\documentclass[11pt]{article}
\usepackage{latexsym}
\usepackage{amsfonts}
\usepackage{mathrsfs,amsthm,amsmath}
\usepackage{color}
\newtheorem{theorem}{Theorem}[section]
\newtheorem{proposition}[theorem]{Proposition}

\newtheorem{lemma}{Lemma}[section]

\newtheorem{remark}{Remark}[section]

\newtheorem{condition}{Condition}[section]
\begin{document}
\title{Asymptotic results for families of power series 
distributions\thanks{All the authors acknowledge the support of 
Indam-GNAMPA (research project \lq\lq Stime asintotiche: principi di 
invarianza e grandi deviazioni\rq\rq). Claudio Macci and Barbara 
Pacchiarotti also acknowledge the support of the MIUR Excellence 
Department Project awarded to the Department of Mathematics, 
University of Rome Tor Vergata (CUP E83C18000100006) and of 
University of Rome Tor Vergata (research program \lq\lq Beyond Borders\rq\rq,
project \lq\lq Asymptotic Methods in Probability\rq\rq (CUP E89C20000680005)).}}
\author{Claudio Macci\thanks{Address: Dipartimento di Matematica,
Universit\`a di Roma Tor Vergata, Via della Ricerca Scientifica,
I-00133 Rome, Italy. e-mail: \texttt{macci@mat.uniroma2.it}}\and
Barbara Pacchiarotti\thanks{Address: Dipartimento di Matematica,
Universit\`a di Roma Tor Vergata, Via della Ricerca Scientifica,
I-00133 Rome, Italy. e-mail: \texttt{pacchiar@mat.uniroma2.it}}\and
Elena Villa\thanks{Address: Dipartimento di Matematica,
Universit\`a degli Studi di Milano, Via Cesare Saldini 50, 
I-20133 Milan, Italy. e-mail: \texttt{elena.villa@unimi.it}}}
%\date{}
\maketitle
\begin{abstract}
In this paper we consider suitable families of power series distributed random variables, and we 
study their asymptotic behavior in the fashion of large (and moderate) deviations. We also present
two examples of fractional counting processes, where the normalizations of the 
involved power series distributions can be expressed in terms of the Prabhakar function. The first
example allows to consider the counting process in \cite{PoganyTomovski}, the second one is inspired 
by a model studied in \cite{GarraOrsingherPolito}.
\end{abstract}

\section{Introduction}
Several discrete distributions in probability concern nonnegative integer valued random variables.
A random variable $X$ has a power series distribution if
$$P(X=k)=\frac{d_k\delta^k}{D(\delta)}\ \mbox{for each integer}\ k\geq 0,$$
where $\delta>0$ is called power parameter, $\{d_k:k\geq 0\}$ is a family of nonnegative numbers 
and the normalization $D(\delta):=\sum_{k\geq 0}d_k\delta^k\in(0,\infty)$ is called the series 
function. Typically analytical properties of the series function $D(\cdot)$ can be related to some
statistical properties of the power series distribution. Moreover the probability generating 
function of a power series distributed random variables $X$ can be easily expressed in terms of the
function $D$; in fact we have
$$\mathbb{E}[u^X]=\frac{D(u\delta)}{D(\delta)}\ \mbox{for all}\ u>0.$$
The reference \cite{Patil} made great advances in the theory of power series distributions. Another 
important contribution was given by the modified power series distributions in \cite{Gupta}, which 
include distributions derived from Lagrangian expansions (see e.g. \cite{ConsulShenton}). 
Other more recent references on these distributions concern some families which contain the geometric
distribution as a particular case: the generalized hypergeometric family, the $q$-series family and 
the Lerch family. Among the references on the Lerch family we recall \cite{GuptaGuptaOngSrivastava}
and \cite{Kemp}; see also \cite{LuoParmarRaina} as a reference on the related Hurwitz-Lerch zeta 
function.

In this paper we consider a family of random variables $\{N(t):t\geq 0\}$, whose univariate 
marginal distributions are expressed in terms of a family of power series distributions 
$\{\mathcal{P}_j:j\geq 0\}$ with power parameter $\delta$; moreover, for all $j\geq 0$, we set 
$\delta:=\delta_j(t)$ for some functions $\{\delta_j(\cdot):j\geq 0\}$. A precise definition 
is given at the beginning of Section \ref{sec:eventually-constant} (some assumptions are needed 
and they are collected in Condition \ref{cond:eventually-constant}) and it is a generalization 
of the basic model with a unique power series distribution, i.e. the case with
$$P(N(t)=k)=\frac{d_k(\delta(t))^k}{D(\delta(t))}\ \mbox{for each integer}\ k\geq 0$$
for some coefficients $\{d_k:k\geq 0\}$, a series function $D(\cdot)$ and a function 
$\delta(\cdot)$.

We recall that, when we deal with the basic model, we have suitable weighted Poisson distributed 
random variables; in fact, for each integer $k\geq 0$, we have
$$P(N(t)=k)=\frac{w(k)\frac{(\lambda\delta(t))^k}{k!}e^{-\lambda\delta(t)}}
{\sum_{j\geq 0}w(j)\frac{(\lambda\delta(t))^j}{j!}e^{-\lambda\delta(t)}},\ \mbox{with}\ w(k)=\frac{k!}{\lambda^k}d_k.$$
This kind of structure was already highlighted in \cite[Section 4]{BeghinMacci2013} for the case
$\delta(t)=t^\nu$ and $d_k=\frac{\lambda^k}{\Gamma(\nu k+1)}$ (for some $\nu\in(0,1]$), and therefore
$w(k)=\frac{k!}{\Gamma(\nu k+1)}$. Weighted Poisson distributions are often related to the concepts
of overdispersion and underdispersion; for some insights on this topic see e.g. \cite{DCPC1998},
\cite{DCPC2005}, the recent paper \cite{CahoyDinardoPolito} and references therein.

Our main results in the present paper concern large (and moderate) deviations as $t\to\infty$ for the 
above mentioned general family $\{N(t):t\geq 0\}$. We remind that the theory of large deviations deals 
with asymptotic computation of small probabilities on an exponential scale (see e.g. \cite{DemboZeitouni} 
as a reference on this topic).

At the best of our knowledge we are not aware of similar results in the literature for power series distributions;
therefore we think that our general results may find applications in several models and may bring to further 
investigations. In this paper we apply the results to different classes of fractional counting processes in the 
literature, where the function $D(\cdot)$ can be expressed in terms of the Prabhakar function (the definition of 
this function will be recalled in Section \ref{sub:ML-Prabhakar-functions}). Namely, in Section \ref{sec:examples} 
we shall present two particular examples. The first one is related to the fractional process in \cite{PoganyTomovski}, 
and it allows us to generalize some large deviation results in the current literature, as discussed at the end of 
Section \ref{sub:basic-model-appl}. The second one is related to a fractional process in \cite{GarraOrsingherPolito}, 
and we discuss a class of cases for which certain conditions on some involved parameters fail.

We also point out that the model in \cite{PoganyTomovski} is a particular case of the family in 
\cite[Section 3.1, eq. (48)]{CahoyDinardoPolito}, where the weights are expressed with a ratio of Gamma functions;
thus, as a possible future work, one might try to investigate a wider class of models defined by suitable generalizations
of the Prabhakar function.

We conclude with the outline of the paper. We start with some preliminaries in Section
\ref{sec:preliminaries}. In Section \ref{sec:eventually-constant} we give a precise
definition of the model, and we prove the results. Finally, in Section \ref{sec:examples},
we apply our results to some examples of fractional counting processes in the literature.

\section{Preliminaries}\label{sec:preliminaries}
In this section we start with some preliminaries on large deviations. Moreover, in view of
the examples presented in Section \ref{sec:examples}, we present some preliminaries on some
special functions.

\subsection{On large deviations}\label{sub:LD}
We start with the definition of large deviation principle (LDP from now on). 
In view of what follows our presentation concerns the case $t\to\infty$; moreover,
for simplicity, we refer to a family of real valued random variables $\{X_t:t>0\}$
defined on the same probability space $(\Omega,\mathcal{F},P)$.

A lower semi-continuous function $I:\mathbb{R}\to[0,\infty]$ is called 
rate function, and it is said to be good if all its level sets
$\{\{x\in\mathbb{R}:I(x)\leq\eta\}:\eta\geq 0\}$ are compact. Then $\{X_t:t>0\}$ 
satisfies the LDP with speed $v_t\to\infty$ and rate function $I$ if
$$\limsup_{t\to\infty}\frac{1}{v_t}\log P(X_t\in C)\leq-\inf_{x\in C}I(x)\ \mbox{for all closed sets}\ C$$
and
$$\liminf_{t\to\infty}\frac{1}{v_t}\log P(X_t\in O)\geq-\inf_{x\in O}I(x)\ \mbox{for all open sets}\ O.$$

We talk about moderate deviations when we have a class of LDPs for families of centered
(or asymptotically centered) random variables which depends on some positive scaling 
factors $\{a(t):t>0\}$ such that
\begin{equation}\label{eq:MD-conditions}
a(t)\to 0\ \mbox{and}\ v_ta(t)\to\infty\ \mbox{as}\ t\to\infty
\end{equation}
and, moreover, all these LDPs (whose speed functions depend on the scaling factors)
are governed by the same quadratic rate function vanishing at zero. We can also say
that, as usually happens, this class of LDPs fills the gap between a convergence to
zero and an asymptotic normality result (see Remark \ref{rem:MD-eventually-constant-fill-the-gap}).

The main large deviation tool used in this paper is the G\"artner Ellis
Theorem (see e.g. \cite[Theorem 2.3.6]{DemboZeitouni}; actually we can refer
to the statement (c) only), and here we recall its statement for real valued random 
variables. In view of this we also recall that a convex function 
$f:\mathbb{R}\to(-\infty,\infty]$ is essentially smooth (see e.g. 
\cite[Definition 2.3.5]{DemboZeitouni}) if the interior of $\mathcal{D}_f:=\{\theta\in\mathbb{R}: 
f(\theta)<\infty\}$ is non-empty, $f$ is differentiable throughout the interior of 
$\mathcal{D}_f$, and $f$ is steep (i.e. $|f^\prime(t)|$ is divergent as $t$ approaches to 
any finite point of the boundary of $\mathcal{D}_f$). In our applications the function $f$ is 
always finite everywhere and differentiable; therefore $f$ is essentially smooth because
the steepness condition holds vacuously.

\begin{theorem}\label{th:GE}
Let $\{X_t:t>0\}$ be a family of real valued random variables defined on the 
same probability space $(\Omega,\mathcal{F},P)$ and let $v_t$ be such that $v_t\to\infty$. 
Moreover assume that, for all $\theta\in\mathbb{R}$, there exists
$$f(\theta):=\lim_{t\to\infty}\frac{1}{v_t}\log\mathbb{E}\left[e^{\theta X_t}\right]$$
as an extended real number; we also assume that the origin $\theta=0$ belongs 
to the interior of the set $\mathcal{D}(f):=\{\theta\in\mathbb{R}:f(\theta)<\infty\}$.
Then, if $f$ is essentially smooth and lower semi-continuous, the family
of random variables $\{X_t/v_t:t>0\}$ satisfies the LDP with speed $v_t$ and good 
rate function $f^*$ defined by $f^*(x):=\sup_{\theta\in\mathbb{R}}\{\theta x-f(\theta)\}$.
\end{theorem}

\subsection{On special functions for some fractional counting processes}\label{sub:ML-Prabhakar-functions}
In this paper, for $\alpha\in(0,1]$ and $\beta,\gamma>0$, we consider the Prabhakar function
$E_{\alpha,\beta}^\gamma(\cdot)$ defined by
$$E_{\alpha,\beta}^\gamma(u):=\sum_{k\geq 0}\frac{u^k(\gamma)_k}{k!\Gamma(\alpha k+\beta)}\ (\mbox{for}\ u\in\mathbb{R}),$$
where
$$(\gamma)_k:=\left\{\begin{array}{ll}
1&\ \mbox{if}\ k=0\\
\gamma(\gamma+1)\cdots(\gamma+k-1)&\ \mbox{if}\ k\geq 1
\end{array}\right.$$
is the rising factorial (Pochhammer symbol). The Prabhakar function is also known as the Mittag-Leffler 
function with three parameters; the Mittag-Leffler function with two parameters concerns the case 
$\gamma=1$, and the classical Mittag-Leffler function concerns the case $\beta=\gamma=1$. Here we are 
interested to the case of a positive argument $u$ and we refer to the asymptotic behavior of 
$E_{\alpha,\beta}^\gamma(\cdot)$ as the argument tends to infinity (see e.g. \cite[page 23]{Giusti-etal}, 
which concerns a result in \cite{Paris} where the argument $z$ of $E_{\alpha,\beta}^\gamma(\cdot)$ is 
complex; obviously we are interested in the case $|\arg(z)|<\frac{\alpha\pi}{2}$). In particular, for 
some $\omega(u)$ such that $\omega(u)\to 0$ as $u\to\infty$, we have 
%$$E_{\alpha,\beta}^\gamma(u)\sim\mathcal{E}_{\alpha,\beta}^\gamma(u)+\mathcal{R}_{\alpha,\beta}^\gamma(u),$$
%where $\sim$ means that the ratio tends to 1 (as $u\to\infty$), $\mathcal{E}_{\alpha,\beta}^\gamma(u)$ 
%is an exponentially large term (as $u\to\infty$) which can be written as follows
\begin{equation}\label{eq:exponentially-large-term}
E_{\alpha,\beta}^\gamma(u)=\frac{1}{\Gamma(\gamma)}e^{u^{1/\alpha}}
u^{\frac{\gamma-\beta}{\alpha}}\frac{1}{\alpha^\gamma}\sum_{k\geq 0}c_ku^{-\frac{k}{\alpha}}(1+\omega(u)),
\end{equation}
where the coefficients $\{c_k:k\geq 0\}$ are obtained by a suitable inverse factorial expansion. Moreover, when 
we present the first application of our results to some fractional counting processes (see Section 
\ref{sub:basic-model-appl}), we shall restrict the attention to the case with a positive integer
$\gamma$; then we refer to \cite[eq. (4.4)]{Giusti-etal}, i.e.
\begin{equation}\label{eq:4.4-Giusti-et-al}
E_{\alpha,\beta}^{\gamma+1}(u)=\frac{1}{\alpha^\gamma \gamma!}\sum_{j=0}^\gamma d_{j,\alpha,\beta}^{(\gamma)}E_{\alpha,\beta-j}^1(u)
\end{equation}
for some coefficients $\{d_{j,\alpha,\beta}^{(\gamma)}:j\in\{0,1,\ldots,\gamma\}\}$ defined by a recursive
expression provided by \cite[eq. (4.6)]{Giusti-etal}, and in particular we have 
$d_{\gamma,\alpha,\beta}^{(\gamma)}=1$. We also recall the following asymptotic formula for the case $\gamma=1$
(see e.g. \cite[eq. (4.4.16)]{GorenfloKilbasMainardiRogosin}), i.e.
\begin{equation}\label{eq:4.4.16-Gorenflo-et-al}
E_{\alpha,\beta}^1(u)=\frac{1}{\alpha}u^{\frac{1-\beta}{\alpha}}e^{u^{1/\alpha}}+O(1/u)\ \mbox{as}\ u\to\infty.
\end{equation}
As we shall see, the Prabhakar function plays an important role in the examples presented in Section 
\ref{sec:examples} based on some fractional counting processes in the literature. We also mention that a 
different use of the Prabhakar function can be found in \cite{GajdaBeghin} for the definition of a new class 
of L\'evy processes (called Prabhakar L\'evy processes).

\section{Model and results}\label{sec:eventually-constant}
We consider a family of power series distributions $\{\mathcal{P}_j:j\geq 0\}$ such that, 
for each $j\geq 0$, $\mathcal{P}_j$ concerns the probability mass function
$$p_j(k):=\frac{d_{k,j}\delta^k}{D_j(\delta)}\ \mbox{for each integer}\ k\geq 0,$$
where $\delta>0$, and $\{d_{k,j}:k\geq 0\}$ is a sequence of nonnegative numbers 
such that
$$D_j(\delta):=\sum_{k\geq 0}d_{k,j}\delta^k\in(0,\infty).$$
Then we consider a family of random variables $\{N(t):t\geq 0\}$ whose probability 
mass functions depend on $\{d_{k,k}:k\geq 0\}$ only; more precisely, assuming that
$$\sum_{j\geq 0}\frac{d_{j,j}\delta^j}{D_j(\delta)}\in(0,\infty)\ \mbox{for all}\ \delta>0,$$
we have
$$P(N(t)=k):=\frac{\frac{d_{k,k}(\delta_k(t))^k}{D_k(\delta_k(t))}}{\sum_{j\geq 0}\frac{d_{j,j}(\delta_j(t))^j}{D_j(\delta_j(t))}}
\ \mbox{for each integer}\ k\geq 0,$$
and $\delta_j(t)\to\infty$ as $t\to\infty$ (for all $j\geq 0$).

\begin{remark}\label{rem:not-proper-stochastic-process}
In the next Section \ref{sec:examples} we show how our results can be applied to some stochastic processes in some references
in the literature. Actually the authors of these references use the term \lq\lq stochastic process\rq\rq\ even if they do not 
specify the joint distributions of the involved random variables. In our results we do not need to define the joint distributions
of the random variables $\{N(t):t\geq 0\}$; indeed we only need to consider their univariate marginal distributions and, for 
this reason, we use the term \lq\lq family of random variables\rq\rq. In a possible future work one could investigate the 
possibility to define all the finite-dimensional distributions in order to prove sample-path large deviation results.	
\end{remark}

In our results some hypotheses are needed, and they are collected in the next Condition
\ref{cond:eventually-constant}.

\begin{condition}\label{cond:eventually-constant}
We consider the following hypotheses.\\
$\mathbf{(B1)}$: There exists $n\geq 0$ such that the elements of both sequences $\{\mathcal{P}_j:j\geq 0\}$ 
and $\{\delta_j(\cdot):j\geq 0\}$ do not depend on $j\geq n$; in particular, for $j\geq n$,
we simply write $d_k$ in place of $d_{k,j}$, $D(\cdot)$ in place of $D_j(\cdot)$, and 
$\delta(\cdot)$ in place of $\delta_j(\cdot)$. Moreover we assume that there exist two functions 
$v:(0,\infty)\to(0,\infty)$ and $\Delta:(0,\infty)\to\mathbb{R}$ such that $v(t)\to\infty$ as $t\to\infty$,
\begin{equation}\label{eq:central}
\lim_{t\to\infty}\frac{1}{v(t)}\log D(ut)=\Delta(u)\ \mbox{for all}\ u>0,
\end{equation}
and $\Delta(\cdot)$ is a differentiable function.\\
$\mathbf{(B2)}$: The set $\{k\geq 0:d_{k,k}>0\}$ is unbounded; thus, if we refer to the case $k\geq n$, the 
set $\{k\geq 0:d_k>0\}$ is unbounded.\\
$\mathbf{(B3)}$: For all $k\in\{0,1,\ldots,n-1\}$ we have:
\begin{equation}\label{eq:limit-ratio}
\lim_{t\to\infty}\frac{\frac{d_{k,k}(\delta_k(t))^k}{D_k(\delta_k(t))}}
{\frac{d_k(\delta(t))^k}{D(\delta(t))}}=0\ \mbox{if}\ d_k>0,
\end{equation}
and $d_{k,k}=0$ if $d_k=0$.
\end{condition}

Firstly we note that the function $\Delta(\cdot)$ is increasing. Moreover, if $n=0$, we have the basic model
with a unique power series distribution and, in particular, $\mathbf{(B3)}$ holds vacuously (because the set 
$\{0,1,\ldots,n-1\}$ in $\mathbf{(B3)}$ is empty).

\begin{remark}\label{rem:finite-sum-extended}
Here we illustrate two consequences of $\mathbf{(B2)}$ in Condition \ref{cond:eventually-constant}. Firstly 
such condition allows to avoid the case $P(0\leq N(t)\leq M)=1$ (for all $t\geq 0$) for 
some $M\in(0,\infty)$; in fact, in such a case, it is easy to check that the results proved below hold with
$\Lambda(\theta)=0$ for all $\theta\in\mathbb{R}$, where $\Lambda$ is the function defined 
in eq. \eqref{eq:Lambda-Lambdastar}. Moreover $\mathbf{(B2)}$ yields
\begin{equation}\label{eq:limit-pmf-k-leq-n-1}
\lim_{t\to\infty}\frac{d_{k,k}(\delta_k(t))^k}{D_k(\delta_k(t))}=0\ \mbox{and}\ \lim_{t\to\infty}\frac{d_k(\delta(t))^k}{D(\delta(t))}=0
\ \mbox{for all}\ k\geq 0;
\end{equation}
in fact, for all $k\geq 0$, there exist $h_1,h_2>k$ such that $d_{h_1,h_1},d_{h_2}>0$; in fact
$$0\leq\frac{d_{k,k}(\delta_k(t))^k}{D_k(\delta_k(t))}\leq\frac{d_{k,k}(\delta_k(t))^k}{d_{h_1,k}(\delta_k(t))^{h_1}}\to 0
\ \mbox{and}\ 0\leq\frac{d_k(\delta(t))^k}{D(\delta(t))}\leq\frac{d_k(\delta(t))^k}{d_{h_2}(\delta(t))^{h_2}}\to 0
\ (\mbox{as}\ t\to\infty).$$
\end{remark}

We also briefly discuss some particular cases concerning the functions $v(\cdot)$ and $\Delta(\cdot)$
in eq. \eqref{eq:central}.

\begin{remark}\label{rem:particular-cases-for-eq:central}
Assume that there exists
\begin{equation}\label{eq:condition-on-v}
\lim_{t\to\infty}\frac{v(ut)}{v(t)}=:\bar{v}(u)\ \mbox{for all}\ u>0
\end{equation}
as a finite limit; then the limit in eq. \eqref{eq:central} can be checked
only for $u=1$, and we have
$$\Delta(u)=\Delta(1)\bar{v}(u)\ \mbox{for all}\ u>0.$$
In particular, if $v(\cdot)$ is a regularly varying function of index 
$\varrho>0$ (see e.g. \cite[Definition A3.1(b)]{EmbrechtsKluppelbergMikosch}),
the limit in eq. \eqref{eq:condition-on-v} holds with $\bar{v}(u)=u^\varrho$.
On the other hand, if $v(\cdot)$ is a slowly varying function (see e.g. 
\cite[Definition A3.1(a)]{EmbrechtsKluppelbergMikosch}), the limit in eq. 
\eqref{eq:condition-on-v} holds with $\bar{v}(u)=1$ (and this case is not 
interesting).
\end{remark}

In view of what follows it is useful to introduce the following notation:
\begin{equation}\label{eq:def-Rn}
R_n(u,t):=\left\{\begin{array}{ll}
0&\ \mbox{if}\ n=0\\
\sum_{k=0}^{n-1}\left(\frac{d_{k,k}(u\delta_k(t))^k}{D_k(\delta_k(t))}-\frac{d_k(u\delta(t))^k}{D(\delta(t))}\right)&\ \mbox{if}\ n\geq 1,
\end{array}\right.
\end{equation}
where $n$ is the value in Condition \ref{cond:eventually-constant}. Then we have
\begin{equation}\label{eq:limit-Rn}
\lim_{t\to\infty}R_n(u,t)=0;
\end{equation}
this is trivial if $n=0$ and, if $n\geq 1$, this is a consequence of eq. 
\eqref{eq:limit-pmf-k-leq-n-1} (for $k\in\{0,1,\ldots,n-1\}$).

We start with the first result.

\begin{proposition}\label{prop:LD-eventually-constant}
Assume that Condition \ref{cond:eventually-constant} holds. Then $\left\{\frac{N(t)}{v(\delta(t))}:t>0\right\}$
satisfies the LDP with speed $v(\delta(t))$ and good rate function $\Lambda^*$ defined by
\begin{equation}\label{eq:Lambda-Lambdastar}
\Lambda^*(x):=\sup_{\theta\in\mathbb{R}}\{\theta x-\Lambda(\theta)\},\ \mbox{where}\ 
\Lambda(\theta):=\Delta(e^\theta)-\Delta(1).
\end{equation}
\end{proposition}
\begin{proof}
We want to apply the G\"artner Ellis Theorem (Theorem \ref{th:GE}). In order to do this
we remark that, for all $\theta\in\mathbb{R}$, we have
\begin{multline*}
\frac{1}{v(\delta(t))}\log\mathbb{E}\left[e^{\theta N(t)}\right]
=\frac{1}{v(\delta(t))}\log\frac{\sum_{k\geq 0}\frac{d_{k,k}(e^\theta\delta_k(t))^k}
{D_k(\delta_k(t))}}{\sum_{j\geq 0}\frac{d_{j,j}(\delta_j(t))^j}{D_j(\delta_j(t))}}\\
=\frac{1}{v(\delta(t))}\log\sum_{k\geq 0}\frac{d_{k,k}(e^\theta\delta_k(t))^k}
{D_k(\delta_k(t))}-\frac{1}{v(\delta(t))}\log\sum_{j\geq 0}\frac{d_{j,j}(\delta_j(t))^j}{D_j(\delta_j(t))}\\
=\frac{1}{v(\delta(t))}\log\left(D(\delta(t))\sum_{k\geq 0}\frac{d_{k,k}(e^\theta\delta_k(t))^k}
{D_k(\delta_k(t))}\right)-\frac{1}{v(\delta(t))}\log\left(D(\delta(t))\sum_{j\geq 0}\frac{d_{j,j}(\delta_j(t))^j}{D_j(\delta_j(t))}\right).
\end{multline*}
So, if we prove that
\begin{equation}\label{eq:central-extended}
\lim_{t\to\infty}\frac{1}{v(\delta(t))}\log\left(D(\delta(t))\sum_{k\geq 0}\frac{d_{k,k}(u\delta_k(t))^k}{D_k(\delta_k(t))}\right)
=\Delta(u)\ \mbox{for all}\ u>0,
\end{equation}
where $\Delta(\cdot)$ is the function in Condition \ref{cond:eventually-constant}, the limit in
eq. \eqref{eq:central-extended} with $u=e^\theta$ and $u=1$ yields
\begin{equation}\label{eq:LD-GETlimit-eventually-constant}
\lim_{t\to\infty}\frac{1}{v(\delta(t))}\log\mathbb{E}\left[e^{\theta N(t)}\right]=\Delta(e^\theta)-\Delta(1)=\Lambda(\theta)
\ \mbox{for all}\ \theta\in\mathbb{R},
\end{equation}
where $\Lambda(\cdot)$ is the function in eq. \eqref{eq:Lambda-Lambdastar}. Then the
desired LDP holds as a straighforward application of Theorem \ref{th:GE}.

So in the remaining part of the proof we show that the limit in eq. 
\eqref{eq:central-extended} holds. This will be done by considering $n\geq 1$; actually, for
$n=0$, we have the same computations and some parts are even simplified. Firstly, if we consider
the function $R_n(u,t)$ defined in eq. \eqref{eq:def-Rn}, for all $u>0$ we have
\begin{equation}\label{eq:initial-step}
\sum_{k\geq 0}\frac{d_{k,k}(u\delta_k(t))^k}{D_k(\delta_k(t))}=\frac{D(u\delta(t))}{D(\delta(t))}+R_n(u,t);
\end{equation}
in fact we have
$$\sum_{k\geq 0}\frac{d_{k,k}(u\delta_k(t))^k}{D_k(\delta_k(t))}
=\sum_{k=0}^{n-1}\frac{d_{k,k}(u\delta_k(t))^k}{D_k(\delta_k(t))}+
\sum_{k\geq n}\frac{d_k(u\delta(t))^k}{D(\delta(t))}\\
=\sum_{k\geq 0}\frac{d_k(u\delta(t))^k}{D(\delta(t))}+R_n(u,t),$$
and we get eq. \eqref{eq:initial-step} by taking into account the definition of
$D(\cdot)$. Then eq. \eqref{eq:initial-step} yields
$$D(\delta(t))\sum_{k\geq 0}\frac{d_{k,k}(u\delta_k(t))^k}{D_k(\delta_k(t))}
=D(u\delta(t))+R_n(u,t)D(\delta(t))=D(u\delta(t))\left(1+R_n(u,t)\frac{D(\delta(t))}{D(u\delta(t))}\right);$$
thus (if we take the logarithms, we divide by $v(\delta(t))$ and we let $t$ go to
infinity) the limit in eq. \eqref{eq:central-extended} holds if we show that
\begin{equation}\label{eq:implies-central-extended}
\lim_{t\to\infty}R_n(u,t)\frac{D(\delta(t))}{D(u\delta(t))}=0\ \mbox{for all}\ u>0.
\end{equation}

So we complete the proof by showing that the limit in eq. \eqref{eq:implies-central-extended} holds.
In fact we have
\begin{multline*}
R_n(u,t)\frac{D(\delta(t))}{D(u\delta(t))}
=\sum_{k=0}^{n-1}\left(\frac{d_{k,k}(u\delta_k(t))^k}{D_k(\delta_k(t))}-\frac{d_k(u\delta(t))^k}{D(\delta(t))}\right)\frac{D(\delta(t))}{D(u\delta(t))}\\
=\sum_{k=0}^{n-1}\left(\frac{\frac{d_{k,k}(u\delta_k(t))^k}{D_k(\delta_k(t))}}{\frac{d_k(u\delta(t))^k}{D(\delta(t))}}
\frac{d_k(u\delta(t))^k}{D(\delta(t))}-\frac{d_k(u\delta(t))^k}{D(\delta(t))}\right)\frac{D(\delta(t))}{D(u\delta(t))}\\
=\sum_{k=0}^{n-1}\left(\frac{\frac{d_{k,k}(\delta_k(t))^k}{D_k(\delta_k(t))}}{\frac{d_k(\delta(t))^k}{D(\delta(t))}}-1\right)
\frac{d_k(u\delta(t))^k}{D(\delta(t))}\frac{D(\delta(t))}{D(u\delta(t))}
=\sum_{k=0}^{n-1}\left(\frac{\frac{d_{k,k}(\delta_k(t))^k}{D_k(\delta_k(t))}}{\frac{d_k(\delta(t))^k}{D(\delta(t))}}-1\right)
\frac{d_k(u\delta(t))^k}{D(u\delta(t))},
\end{multline*}
and the desired limit in eq. \eqref{eq:implies-central-extended} holds by the limit in eq. 
\eqref{eq:limit-ratio}, and by the second limit in eq. \eqref{eq:limit-pmf-k-leq-n-1} (here 
we have $u\delta(t)$ instead of $\delta(t)$ and that limit still holds).
\end{proof}

Now we study moderate deviations. More precisely we prove a class of LDPs which depends 
on any possible choice of positive numbers $\{a(t):t>0\}$ such that \eqref{eq:MD-conditions}
holds with $v_t=v(\delta(t))$, which is the speed in Proposition \ref{prop:LD-eventually-constant}.
We remark that $\Lambda^{\prime\prime}(0)$ that appears below (Proposition \ref{prop:MD-eventually-constant}
and Remark \ref{rem:MD-eventually-constant-fill-the-gap}) cannot be negative; in fact, as we have seen in 
the proof of Proposition \ref{prop:LD-eventually-constant}, the function $\Lambda$ is the pointwise limit of 
logarithms of moment generating functions, which are convex functions (see e.g. 
\cite[Lemma 2.2.5(a)]{DemboZeitouni}).

\begin{proposition}\label{prop:MD-eventually-constant}
Assume that Condition \ref{cond:eventually-constant} holds and, if we refer to the function $\Delta(\cdot)$ 
in that condition, let $\Lambda(\cdot)$ be the function in eq. \eqref{eq:Lambda-Lambdastar}. Assume that there 
exists $\Delta^{\prime\prime}(1)$, and therefore there exists $\Lambda^{\prime\prime}(0)$. Moreover assume that,
for $D(\cdot)$, $\delta(\cdot)$ and $v(\cdot)$ in Condition \ref{cond:eventually-constant} (and for $\Lambda(\cdot)$
in eq. \eqref{eq:Lambda-Lambdastar}), the following conditions hold:
\begin{equation}\label{eq:MD-basic-hyp1}
\mbox{if}\ u(t)\to 1\ \mbox{as}\ t\to\infty,\ \mbox{then}\ 
H_1(t):=\log\frac{D\left(u(t)\delta(t)\right)}{D(\delta(t))}-v(\delta(t))(\Delta(u(t))-\Delta(1))\ \mbox{is bounded};
\end{equation}
\begin{equation}\label{eq:MD-basic-hyp2}
H_2(t):=\sqrt{v(\delta(t))}\left(\Lambda^\prime(0)-\frac{\delta(t)D^\prime(\delta(t))}{v(\delta(t))D(\delta(t))}\right)\ \mbox{is bounded};
\end{equation}
\begin{equation}\label{eq:extra-for-MD-eventually-constant}
H_3(t):=\frac{1}{\sqrt{v(\delta(t))}}\left(\frac{\delta(t)D^\prime(\delta(t))}{D(\delta(t))}-\mathbb{E}[N(t)]\right)\ \mbox{is bounded}.
\end{equation}
Then, for every choice of $\{a(t):t>0\}$ such that eq. \eqref{eq:MD-conditions} holds with $v_t=v(\delta(t))$, the 
random variables 
$\left\{\frac{N(t)-\mathbb{E}[N(t)]}{v(\delta(t))}\sqrt{v(\delta(t))a(t)}:t>0\right\}$
satisfies the LDP with speed $1/a(t)$ and good rate function $\tilde{\Lambda}^*$ defined by
\begin{equation}\label{eq:tilde-Lambda-star}
\tilde{\Lambda}^*(x):=\left\{\begin{array}{ll}
\frac{x^2}{2\Lambda^{\prime\prime}(0)}&\ \mbox{for}\ \Lambda^{\prime\prime}(0)>0\\
\left\{\begin{array}{ll}
0&\ \mbox{if}\ x=0\\
\infty&\ \mbox{if}\ x\neq 0
\end{array}\right.&\ \mbox{for}\ \Lambda^{\prime\prime}(0)=0.
\end{array}\right.
\end{equation}
\end{proposition}
\begin{proof}
We want to apply the G\"artner Ellis Theorem (Theorem \ref{th:GE}). So in what follows we 
show that
\begin{equation}\label{eq:def-Lambda-tilde}
\lim_{t\to\infty}\frac{1}{1/a(t)}\log\mathbb{E}\left[e^{\frac{\theta}{a(t)}
\frac{N(t)-\mathbb{E}[N(t)]}{v(\delta(t))}\sqrt{v(\delta(t))a(t)}}\right]
=\frac{\theta^2}{2}\Lambda^{\prime\prime}(0)=:\tilde{\Lambda}(\theta)\ \mbox{for all}\ \theta\in\mathbb{R};
\end{equation}
in fact it is easy to check that
$$\tilde{\Lambda}^*(x):=\sup_{\theta\in\mathbb{R}}\{\theta x-\tilde{\Lambda}(\theta)\}$$
coincides with the rate function in the statement of the proposition.

Firstly we observe that
\begin{multline*}
\Lambda_t(\theta):=\frac{1}{1/a(t)}\log\mathbb{E}\left[e^{\frac{\theta}{a(t)}\frac{N(t)-\mathbb{E}[N(t)]}
{v(\delta(t))}\sqrt{v(\delta(t))a(t)}}\right]\\
=a(t)\left(\log\frac{\sum_{k\geq 0}\frac{d_{k,k}\left(e^{\frac{\theta}{\sqrt{v(\delta(t))a(t)}}}\delta_k(t)\right)^k}{D_k(\delta_k(t))}}{\sum_{j\geq 0}\frac{d_{j,j}(\delta_j(t))^j}{D_j(\delta_j(t))}}-\frac{\theta}{\sqrt{v(\delta(t))a(t)}}\mathbb{E}[N(t)]\right).
\end{multline*}
Moreover we set again $u(t):=e^{\frac{\theta}{\sqrt{v(\delta(t))a(t)}}}$; in fact, by
eq. \eqref{eq:MD-conditions} with $v_t=v(\delta(t))$, we have $u(t)\to 1$ because
$\sqrt{v(\delta(t))a(t)}\to\infty$. Then we can check that
$$\Lambda_t(\theta)=A_1(t)+A_2(t)+A_3(t),$$
where
\begin{multline*}
A_1(t):=a(t)\left(\log\frac{\sum_{k\geq 0}\frac{d_{k,k}\left(e^{\frac{\theta}{\sqrt{v(\delta(t))a(t)}}}\delta_k(t)\right)^k}{D_k(\delta_k(t))}}
{\sum_{j\geq 0}\frac{d_{j,j}(\delta_j(t))^j}{D_j(\delta_j(t))}}
-v(\delta(t))\Lambda\left(\frac{\theta}{\sqrt{v(\delta(t))a(t)}}\right)\right)\\
=a(t)\left(\log\frac{\frac{D(u(t)\delta(t))}{D(\delta(t))}+R_n(u(t),t)}
{1+R_n(1,t)}-v(\delta(t))\Lambda\left(\frac{\theta}{\sqrt{v(\delta(t))a(t)}}\right)\right)
\end{multline*}
(in the last equality we take into account eq. \eqref{eq:initial-step} with $u=u(t)$ and $u=1$),
$$A_2(t):=v(\delta(t))a(t)\left(\Lambda\left(\frac{\theta}{\sqrt{v(\delta(t))a(t)}}\right)
-\frac{\theta}{v(\delta(t))\sqrt{v(\delta(t))a(t)}}\frac{\delta(t)D^\prime(\delta(t))}{D(\delta(t))}\right),$$
and
$$A_3(t):=a(t)\frac{\theta}{\sqrt{v(\delta(t))a(t)}}\left(\frac{\delta(t)D^\prime(\delta(t))}{D(\delta(t))}-\mathbb{E}[N(t)]\right).$$
So, if we refer to the function $\tilde{\Lambda}(\cdot)$ in eq. \eqref{eq:def-Lambda-tilde}, we 
complete the proof if we show that (for all $\theta\in\mathbb{R}$)
\begin{equation}\label{eq:three-limits}
\lim_{t\to\infty}A_1(t)=0,\ \lim_{t\to\infty}A_2(t)=\tilde{\Lambda}(\theta),\ \lim_{t\to\infty}A_3(t)=0.
\end{equation}

We start by considering $H_1(t)$ in eq. \eqref{eq:MD-basic-hyp1}, and we have
$$H_1(t)=\log\frac{D(u(t)\delta(t))}{D(\delta(t))}-v(\delta(t))\Lambda\left(\frac{\theta}{\sqrt{v(\delta(t))a(t)}}\right)$$
by the definition of the function $\Lambda(\cdot)$ in eq. \eqref{eq:Lambda-Lambdastar}
and by $u(t)=e^{\frac{\theta}{\sqrt{v(\delta(t))a(t)}}}$. Then we can easily check that
\begin{multline*}
A_1(t)=a(t)H_1(t)+a(t)\left(\log\frac{\frac{D(u(t)\delta(t))}{D(\delta(t))}+R_n(u(t),t)}
{1+R_n(1,t)}-\log\frac{D(u(t)\delta(t))}{D(\delta(t))}\right)\\
=a(t)H_1(t)+a(t)\log\left(1+R_n(u(t),t)\frac{D(\delta(t))}{D(u(t)\delta(t))}\right)-a(t)\log(1+R_n(1,t)),
\end{multline*}
where, since $a(t)\to 0$, $a(t)H_1(t)\to 0$ by eq. \eqref{eq:MD-basic-hyp1},
and $a(t)\log(1+R_n(1,t))\to 0$ by eq. \eqref{eq:limit-Rn} with $u=1$. Moreover we have
$$\lim_{t\to\infty}R_n(u(t),t)\frac{D(\delta(t))}{D(u(t)\delta(t))}=0;$$
in fact this is trivial if $n=0$ and, if $n\geq 1$, we have 
\begin{multline*}
0\leq|R_n(u(t),t)|\frac{D(\delta(t))}{D(u(t)\delta(t))}
=\sum_{k=0}^{n-1}\left|\frac{d_{k,k}(u(t)\delta_k(t))^k}{D_k(\delta_k(t))}-\frac{d_k(u(t)\delta(t))^k}{D(\delta(t))}\right|\frac{D(\delta(t))}{D(u(t)\delta(t))}\\
=\sum_{k=0}^{n-1}\left|\frac{\frac{d_{k,k}(u(t)\delta_k(t))^k}{D_k(\delta_k(t))}}{\frac{d_k(u(t)\delta(t))^k}{D(\delta(t))}}-1\right|
\frac{d_k(u(t)\delta(t))^k}{D(u(t)\delta(t))}
=\sum_{k=0}^{n-1}\left|\frac{\frac{d_{k,k}(\delta_k(t))^k}{D_k(\delta_k(t))}}{\frac{d_k(\delta(t))^k}{D(\delta(t))}}-1\right|
\frac{d_k(u(t)\delta(t))^k}{D(u(t)\delta(t))}
\end{multline*}
and, since $u(t)\to 1$, the last expression tends to zero by the limit in eq. \eqref{eq:limit-ratio},
and by the second limit in eq. \eqref{eq:limit-pmf-k-leq-n-1}. Then the first limit in eq.
\eqref{eq:three-limits} is checked.

Now we consider the Taylor formula for $\Lambda(\cdot)$, and we have
$$\Lambda(\eta)=\underbrace{\Lambda(0)}_{=0}+\Lambda^\prime(0)\eta+\frac{\Lambda^{\prime\prime}(0)}{2}\eta^2+o(\eta^2)$$
where $\frac{o(\eta^2)}{\eta^2}\to 0$ as $\eta\to 0$. Then
\begin{multline*}
A_2(t)=v(\delta(t))a(t)\left(\Lambda\left(\frac{\theta}{\sqrt{v(\delta(t))a(t)}}\right)
-\frac{\theta}{\sqrt{v(\delta(t))a(t)}}\frac{\delta(t)D^\prime(\delta(t))}{v(\delta(t))D(\delta(t))}\right)\\
=v(\delta(t))a(t)\left(\left(\Lambda^\prime(0)-\frac{\delta(t)D^\prime(\delta(t))}{v(\delta(t))D(\delta(t))}\right)
\frac{\theta}{\sqrt{v(\delta(t))a(t)}}+\frac{\Lambda^{\prime\prime}(0)}{2}\frac{\theta^2}{v(\delta(t))a(t)}
+o\left(\frac{1}{v(\delta(t))a(t)}\right)\right)\\
=\sqrt{a(t)}\theta H_2(t)+\frac{\Lambda^{\prime\prime}(0)}{2}\theta^2+v(\delta(t))a(t)o\left(\frac{1}{v(\delta(t))a(t)}\right),
\end{multline*}
and the second limit in eq. \eqref{eq:three-limits} holds by eq. \eqref{eq:MD-basic-hyp2} and $a(t)\to 0$, and also by 
$v(\delta(t))a(t)\to\infty$.

Finally we have
$$A_3(t)=\sqrt{a(t)}\theta H_3(t)$$
and the third limit in eq. \eqref{eq:three-limits} holds by eq. \eqref{eq:extra-for-MD-eventually-constant} and $a(t)\to 0$.
\end{proof}

We conclude with some consequences of Proposition \ref{prop:MD-eventually-constant}, 
which are typical features of moderate deviations.

\begin{remark}\label{rem:MD-eventually-constant-fill-the-gap}
The class of LDPs in Proposition \ref{prop:MD-eventually-constant} fill the gap between two
following asymptotic regimes.
\begin{enumerate}
\item The weak convergence of 
$\left\{\frac{N(t)-\mathbb{E}[N(t)]}{\sqrt{v(\delta(t))}}:t>0\right\}$
to the centered Normal distribution with variance $\Lambda^{\prime\prime}(0)$ (in fact the 
proof of Proposition \ref{prop:MD-eventually-constant} still works if $a(t)=1$ and, in such a case, the 
first condition in eq. \eqref{eq:MD-conditions} fails).
\item The convergence of $\left\{\frac{N(t)-\mathbb{E}[N(t)]}{v(\delta(t))}:t>0\right\}$
to zero (in probability) which corresponds to the case $a(t)=\frac{1}{v(\delta(t))}$ 
(in such a case the second condition in eq. \eqref{eq:MD-conditions}, with $v_t=v(\delta(t))$, fails).
\end{enumerate}
Actually in the second case we have in mind cases in which the limit
\begin{equation}\label{eq:limit-for-mean}
\lim_{t\to\infty}\frac{\mathbb{E}[N(t)]}{v(\delta(t))}=\Lambda^\prime(0)
\end{equation}
holds. To better explain this fact we remark that, if the limit in eq. \eqref{eq:limit-for-mean} holds,
then we have
$$\lim_{t\to\infty}\frac{1}{v(\delta(t))}\log\mathbb{E}\left[e^{\theta(N(t)-\mathbb{E}[N(t)])}\right]
=\lim_{t\to\infty}\frac{1}{v(\delta(t))}\log\mathbb{E}\left[e^{\theta N(t)}\right]-\theta\frac{\mathbb{E}[N(t)]}{v(\delta(t))}
=\Lambda(\theta)-\theta\Lambda^\prime(0)$$
for all $\theta\in\mathbb{R}$ (here we take into account the limit in eq. \eqref{eq:LD-GETlimit-eventually-constant}); 
then, if we apply the G\"artner Ellis Theorem (Theorem \ref{th:GE}), the family of random variables 
$\left\{\frac{N(t)-\mathbb{E}[N(t)]}{v(\delta(t))}:t>0\right\}$
satisfies the LDP with speed $v(\delta(t))$ and good rate function $J$ defined by
$$J(y):=\sup_{\theta\in\mathbb{R}}\{\theta y-(\Lambda(\theta)-\theta\Lambda^\prime(0))\}=\Lambda^*(y+\Lambda^\prime(0)),$$
and the rate function $J$ uniquely vanishes at $y=0$ (because $\Lambda^*(x)$ uniquely vanishes at 
$x=\Lambda^\prime(0)$).
\end{remark}

\section{Application of results to some fractional counting processes}\label{sec:examples}
In this section we present two examples of applications of our results to some fractional
counting processes in the literature; so we refer to the content of Section \ref{sub:ML-Prabhakar-functions}.
The first example (in Section \ref{sub:basic-model-appl}) concerns the basic model,
i.e. the case $n=0$; the second example (in Section \ref{sub:eventually-constant-appl}) depends on
two sequences of parameters $\{\alpha_j:j\geq 0\}$ and $\{\tilde{\alpha}_j:j\geq 0\}$ satisfying
suitable conditions. So, in Section \ref{sub:counterexample}, we discuss a class of cases for which
such conditions fail, and we cannot refer to a straightforward application of our results because 
the hypotheses of the G\"artner Ellis Theorem (Theorem \ref{th:GE}) fail.

\subsection{An example related to the basic model}\label{sub:basic-model-appl}
A reference for this example is \cite{PoganyTomovski}; some other connections with literature are
presented below in the last paragraph of this section. In particular it is a case with $n=0$. 
For $\beta,\gamma,\lambda>0$ and $\alpha\in(0,1]$, we set
$$d_k:=\frac{\lambda^k(\gamma)_k}{k!\Gamma(\alpha k+\beta)},$$
where $(\gamma)_k$ is the rising factorial; therefore we get
$$D(u)=E_{\alpha,\beta}^\gamma(\lambda u),$$
where $E_{\alpha,\beta}^\gamma(\cdot)$ is the Prabhakar function.

We start with a discussion on Condition \ref{cond:eventually-constant}. Moreover we discuss eqs. 
\eqref{eq:MD-basic-hyp1}, \eqref{eq:MD-basic-hyp2} and \eqref{eq:extra-for-MD-eventually-constant} in 
Proposition \ref{prop:MD-eventually-constant}; in this case we assume that $\gamma$ is a positive integer.

\paragraph{Discussion on Condition \ref{cond:eventually-constant}.}
We start noting $\mathbf{(B2)}$ and $\mathbf{(B3)}$ trivially holds. Moreover, as far as $\mathbf{(B1)}$ 
is concerned, we have
\begin{equation}\label{eq:Lambda-for-PoganyTomovski}
v(\delta(t)):=(\delta(t))^{1/\alpha}\ \mbox{and}\ \Delta(u):=(\lambda u)^{1/\alpha},
\ \mbox{and therefore we have}\ \Lambda(\theta)=\lambda^{1/\alpha}(e^{\theta/\alpha}-1)
\end{equation}
(we refer to eq. \eqref{eq:exponentially-large-term} for the limit in eq. \eqref{eq:central}).
Note that the function $v(\cdot)$ in eq. \eqref{eq:Lambda-for-PoganyTomovski} is regularly 
varying with index $\varrho=\frac{1}{\alpha}$; in fact (see Remark \ref{rem:particular-cases-for-eq:central})
we have $\Delta(u)=\Delta(1)\bar{v}(u)$ with $\bar{v}(u)=u^{1/\alpha}$ and $\Delta(1)=\lambda^{1/\alpha}$.

\paragraph{Discussion on eqs. \eqref{eq:MD-basic-hyp1}, \eqref{eq:MD-basic-hyp2} and
\eqref{eq:extra-for-MD-eventually-constant} in Proposition \ref{prop:MD-eventually-constant} 
(when $\gamma$ is a positive integer).}
In view of what follows we remark that, by eqs. \eqref{eq:4.4-Giusti-et-al}-\eqref{eq:4.4.16-Gorenflo-et-al}
and $d_{\gamma,\alpha,\beta}^{(\gamma)}=1$, we have
\begin{multline}
E_{\alpha,\beta}^{\gamma+1}(u)=\frac{1}{\alpha^\gamma \gamma!}
\left(E_{\alpha,\beta-\gamma}^1(u)+\sum_{j=0}^{\gamma-1}d_{j,\alpha,\beta}^{(\gamma)}E_{\alpha,\beta-j}^1(u)\right)\\
=\frac{e^{u^{1/\alpha}}}{\alpha^{\gamma+1}\gamma!}
\left(u^{\frac{\gamma+1-\beta}{\alpha}}+\sum_{j=0}^{\gamma-1}d_{j,\alpha,\beta}^{(\gamma)}u^{\frac{j+1-\beta}{\alpha}}+O(e^{-u^{1/\alpha}}/u)\right)\\
=\frac{e^{u^{1/\alpha}}u^{\frac{\gamma+1-\beta}{\alpha}}}{\alpha^{\gamma+1}\gamma!}
\left(1+\sum_{j=0}^{\gamma-1}d_{j,\alpha,\beta}^{(\gamma)}u^{\frac{j-\gamma}{\alpha}}+o(u^{-1/\alpha})\right)
=\frac{e^{u^{1/\alpha}}u^{\frac{\gamma+1-\beta}{\alpha}}}{\alpha^{\gamma+1}\gamma!}\left(1+O(u^{-1/\alpha})\right).\label{eq:Elena}
\end{multline}

We start with eq. \eqref{eq:MD-basic-hyp1}. We take $u(t)\to 1$ as $t\to\infty$ and, by eq. 
\eqref{eq:Elena}, we have
\begin{multline*}
H_1(t)=\log\frac{E_{\alpha,\beta}^\gamma\left(\lambda u(t)\delta(t)\right)}
{E_{\alpha,\beta}^\gamma(\lambda\delta(t))}-(\delta(t))^{1/\alpha}(\lambda^{1/\alpha}(u(t))^{1/\alpha}-\lambda^{1/\alpha})\\
=\log\frac{e^{(\lambda u(t)\delta(t))^{1/\alpha}}(\lambda u(t)\delta(t))^{\frac{\gamma-\beta}{\alpha}}\left(1+O((u(t)\delta(t))^{-1/\alpha})\right)}
{e^{(\lambda\delta(t))^{1/\alpha}}(\lambda\delta(t))^{\frac{\gamma-\beta}{\alpha}}\left(1+O((\delta(t))^{-1/\alpha})\right)}
-(\lambda\delta(t))^{1/\alpha}((u(t))^{1/\alpha}-1)\\
\frac{\gamma-\beta}{\alpha}\log u(t)+\log\left(1+O((u(t)\delta(t))^{-1/\alpha})\right)
-\log\left(1+O((\delta(t))^{-1/\alpha})\right)\to 0\ (\mbox{as}\ t\to\infty).
\end{multline*}
Thus $H_1(t)$ is bounded and eq. \eqref{eq:MD-basic-hyp1} holds.

Now we consider eq. \eqref{eq:MD-basic-hyp2}. We recall that
$$D^\prime(u)=\lambda\frac{d}{du}E_{\alpha,\beta}^\gamma(\lambda u)=\lambda\gamma E_{\alpha,\alpha+\beta}^{\gamma+1}(\lambda u)$$
(see e.g. \cite[eq. (1.9.5) with $n=1$]{KilbasSrivastavaTrujillo}). Then, since 
$\Lambda^\prime(0)=\frac{\lambda^{1/\alpha}}{\alpha}$, again by eq. \eqref{eq:Elena} we get
\begin{multline*}
H_2(t)=\sqrt{(\delta(t))^{1/\alpha}}\left(\Lambda^\prime(0)-
\frac{\delta(t)\lambda\gamma E_{\alpha,\alpha+\beta}^{\gamma+1}(\lambda\delta(t))}
{(\delta(t))^{1/\alpha}E_{\alpha,\beta}^\gamma(\lambda\delta(t))}\right)\\
=(\delta(t))^{1/(2\alpha)}\left(\frac{\lambda^{1/\alpha}}{\alpha}-
\frac{(\delta(t))^{1-1/\alpha}\lambda\gamma\frac{e^{(\lambda\delta(t))^{1/\alpha}}(\lambda\delta(t))^{\frac{\gamma+1-\alpha-\beta}{\alpha}}}
{\alpha^{\gamma+1}\gamma!}\left(1+O((\delta(t))^{-1/\alpha})\right)}
{\frac{e^{(\lambda\delta(t))^{1/\alpha}}(\lambda\delta(t))^{\frac{\gamma-\beta}{\alpha}}}{\alpha^\gamma(\gamma-1)!}
\left(1+O((\delta(t))^{-1/\alpha})\right)}\right)\\
=(\delta(t))^{1/(2\alpha)}\left(\frac{\lambda^{1/\alpha}}{\alpha}-\frac{\lambda^{1/\alpha}}{\alpha}
\left(\frac{1+O((\delta(t))^{-1/\alpha})}{1+O((\delta(t))^{-1/\alpha})}\right)\right)\\
=\frac{\lambda^{1/\alpha}}{\alpha}(\delta(t))^{1/(2\alpha)}\frac{O((\delta(t))^{-1/\alpha})}{1+O((\delta(t))^{-1/\alpha})}
=\frac{\lambda^{1/\alpha}}{\alpha}\frac{O((\delta(t))^{-1/(2\alpha)})}{1+O((\delta(t))^{-1/\alpha})}\to 0\ (\mbox{as}\ t\to\infty).
\end{multline*}
Thus $H_2(t)$ is bounded and eq. \eqref{eq:MD-basic-hyp2} holds.

\begin{remark}\label{rem:limit-for-mean}
We have just shown that $H_2(t)\to 0$ as $t\to\infty$; then we can immediately check the limit in 
eq. \eqref{eq:limit-for-mean} in Remark \ref{rem:MD-eventually-constant-fill-the-gap} noting that
$$H_2(t)=\sqrt{v(\delta(t))}\left(\Lambda^\prime(0)-\frac{\mathbb{E}[N(t)]}{v(\delta(t))}\right)$$
and $v(\delta(t))\to\infty$.
\end{remark}

We conclude with eq. \eqref{eq:extra-for-MD-eventually-constant} which can be immediately checked;
in fact we have $H_3(t)=0$ because $n=0$, and therefore $H_3(t)$ is bounded.

\paragraph{Connections with the literature.}
If we set $\beta=\gamma=1$ and $\delta(t)=t^\alpha$, we recover the case in \cite[Section 4]{BeghinMacci2013},
and therefore the case in \cite{BeghinMacci2017} with $m=1$. Moreover the function $\Lambda$ in eq. 
\eqref{eq:Lambda-for-PoganyTomovski} coincides with the function $\Lambda_{\alpha,\lambda}(\theta)$ in the proof
of Proposition 4.1 in \cite{BeghinMacci2013} and with the function $\Lambda(\theta)$ in \cite[eq. (7)]{BeghinMacci2017}
specialized to the case $m=1$ (in both cases the parameter $\nu$ in \cite{BeghinMacci2013} and \cite{BeghinMacci2017} 
coincides with $\alpha$ here). In particular we recover the case of Proposition 2 in \cite{BeghinMacci2017} with $m=1$ 
by applying Proposition \ref{prop:MD-eventually-constant} to the example in this section with $\beta=\gamma=1$ and 
$\delta(t)=t^\alpha$.

We also note that, by eq. \eqref{eq:Lambda-for-PoganyTomovski}, we get
$$\Lambda^\prime(0)=\frac{\lambda^{1/\alpha}}{\alpha}\ \mbox{and}\ \Lambda^{\prime\prime}(0)=\frac{\lambda^{1/\alpha}}{\alpha^2}.$$
In particular the equality $\Lambda^{\prime\prime}(0)=\frac{\lambda^{1/\alpha}}{\alpha^2}$ 
coincides with the equality in \cite[eq. (9)]{BeghinMacci2017} for the matrix $C$ 
specialized to the case $m=1$ (and therefore the matrix reduces to a number); in fact,
if we consider $\alpha$ in place of the parameter $\nu$ in \cite{BeghinMacci2017}, we
get
$$c^{(\alpha)}:=\frac{1}{\alpha}\left(\frac{1}{\alpha}-1\right)\lambda^{1/\alpha}+\frac{1}{\alpha}\lambda^{1/\alpha}=\frac{\lambda^{1/\alpha}}{\alpha^2}.$$

\subsection{An example with eventually constant parameters}\label{sub:eventually-constant-appl}
A reference for this example is \cite{GarraOrsingherPolito}; more precisely we refer to
the definition in (3.5) therein. Here we assume that $n\geq 1$; in fact, if $n=0$,
we have a particular case of Example 1. For $\lambda>0$ and for some $\{\alpha_j:j\geq 0\}$ 
with $\alpha_j\in(0,1]$ for all $j\geq 0$, we set
$$d_{k,j}:=\frac{\lambda^k}{\Gamma(\alpha_jk+1)}\ (\mbox{for all}\ k,j\geq 0);$$
therefore we get
$$D_j(u)=E_{\alpha_j,1}^1(\lambda u),$$
where $E_{\alpha_j,1}^1(\cdot)$ is the Mittag-Leffler function (with $\alpha=\alpha_j$).

As far as the functions $\{\delta_j(\cdot):j\geq 0\}$ are concerned, here we consider
the case
$$\delta_j(t):=t^{\tilde{\alpha}_j}$$
for some $\tilde{\alpha}_j\in(0,1]$ (for all $j\geq 0$). Note that the parameters 
$\{\tilde{\alpha}_j:j\geq 0\}$ allow to have a generalization of the case in 
\cite[eq. (3.5)]{GarraOrsingherPolito}, which can be recovered by setting $\tilde{\alpha}_j=1$ 
(for all $j\geq 0$).

\subsubsection{On the conditions in Section \ref{sec:eventually-constant}}\label{sub:discussions}
We start with a discussion on Condition \ref{cond:eventually-constant}. Moreover we discuss eqs. 
\eqref{eq:MD-basic-hyp1}, \eqref{eq:MD-basic-hyp2} and \eqref{eq:extra-for-MD-eventually-constant} in
Proposition \ref{prop:MD-eventually-constant}. In particular, as far as Condition \ref{cond:eventually-constant}
is concerned, we present sufficient conditions on the parameters $\{\alpha_j:j\geq 0\}$ and 
$\{\tilde{\alpha}_j:j\geq 0\}$ in order to have $\mathbf{(B3)}$; moreover, in order to explain what 
can happen when these sufficient conditions fail, a class of cases is studied in detail in the next 
Section \ref{sub:counterexample}.

\paragraph{Discussion on Condition \ref{cond:eventually-constant}.}
We start with $\mathbf{(B1)}$. It is easy to check that we have to consider the following 
restrictions on the parameters that do not appear in \cite{GarraOrsingherPolito}: there 
exist $n\geq 1$ and $\tilde{\alpha},\alpha\in(0,1]$ such that
$$\tilde{\alpha}_j=\tilde{\alpha}\ \mbox{and}\ \alpha_j=\alpha\ \mbox{for all}\ j\geq n.$$
Thus we have
$$\delta(t)=t^{\tilde{\alpha}}$$
and, for $j\geq n$, we can refer to the application to fractional counting processes in Section
\ref{sub:basic-model-appl} with $\beta=\gamma=1$; thus we set
$$d_k:=\frac{\lambda^k}{\Gamma(\alpha k+1)}\ (\mbox{for all}\ k\geq 0),$$
and we have
$$D(u):=E_{\alpha,1}^1(\lambda u).$$
Then, if we refer the statement above with eq. \eqref{eq:Lambda-for-PoganyTomovski} (with 
$\beta=\gamma=1$), we can say that $\mathbf{(B1)}$ holds with 
$$v(t)=t^{1/\alpha}\ \mbox{and}\ \Delta(u)=(\lambda u)^{1/\alpha};$$
thus, in particular, we have
$$v(\delta(t))=t^{\tilde{\alpha}/\alpha}.$$

Condition $\mathbf{(B2)}$ trivially holds because all the coefficients $\{d_{k,j}:k,j\geq 0\}$
are positive. We also note that the limits in eq. \eqref{eq:limit-pmf-k-leq-n-1} hold; in fact
(see Remark \ref{rem:finite-sum-extended}) we have
\begin{equation}\label{eq:limit1-example-pmf-k-leq-n-1}
\frac{d_{k,k}(\delta_k(t))^k}{D_k(\delta_k(t))}
=\frac{\lambda^k}{\Gamma(\alpha_kk+1)}\frac{(t^{\tilde{\alpha}_k})^k}{E_{\alpha_k,1}^1(\lambda t^{\tilde{\alpha}_k})}\to 0
\ (\mbox{as}\ t\to\infty)
\end{equation}
and
\begin{equation}\label{eq:limit2-example-pmf-k-leq-n-1}
\frac{d_k(\delta(t))^k}{D(\delta(t))}
=\frac{\lambda^k}{\Gamma(\alpha k+1)}\frac{(t^{\tilde{\alpha}})^k}{E_{\alpha,1}^1(\lambda t^{\tilde{\alpha}})}\to 0
\ (\mbox{as}\ t\to\infty),
\end{equation}
where the limits hold by eq. \eqref{eq:4.4.16-Gorenflo-et-al} with $u=\lambda t^{\tilde{\alpha}_k}$
and $u=\lambda t^{\tilde{\alpha}}$.

Finally we discuss $\mathbf{(B3)}$. We trivially have $d_k>0$ and, moreover,
$$\frac{\frac{d_{k,k}(\delta_k(t))^k}{D_k(\delta_k(t))}}{\frac{d_k(\delta(t))^k}{D(\delta(t))}}=
\frac{\frac{\lambda^k}{\Gamma(\alpha_kk+1)}\frac{(t^{\tilde{\alpha}_k})^k}{E_{\alpha_k,1}^1(\lambda t^{\tilde{\alpha}_k})}}
{\frac{\lambda^k}{\Gamma(\alpha k+1)}\frac{(t^{\tilde{\alpha}})^k}{E_{\alpha,1}^1(\lambda t^{\tilde{\alpha}})}}
=\frac{\Gamma(\alpha k+1)}{\Gamma(\alpha_kk+1)}t^{(\tilde{\alpha}_k-\tilde{\alpha})k}
\frac{E_{\alpha,1}^1(\lambda t^{\tilde{\alpha}})}{E_{\alpha_k,1}^1(\lambda t^{\tilde{\alpha}_k})};$$
thus, by taking into account again eq. \eqref{eq:4.4.16-Gorenflo-et-al} with 
$u=\lambda t^{\tilde{\alpha}_k}$ and $u=\lambda t^{\tilde{\alpha}}$, the limit in eq. 
\eqref{eq:limit-ratio} holds if, for all $k\in\{0,1,\ldots,n-1\}$, we have
\begin{equation}\label{eq:*}
\frac{\tilde{\alpha}}{\alpha}-\frac{\tilde{\alpha}_k}{\alpha_k}<0
\end{equation}
or
$$\frac{\tilde{\alpha}}{\alpha}-\frac{\tilde{\alpha}_k}{\alpha_k}=0\ \mbox{and}\ \tilde{\alpha}_k-\tilde{\alpha}<0.$$

\paragraph{Discussion on eqs. \eqref{eq:MD-basic-hyp1}, \eqref{eq:MD-basic-hyp2} and
\eqref{eq:extra-for-MD-eventually-constant} in Proposition \ref{prop:MD-eventually-constant}.}
For eqs. \eqref{eq:MD-basic-hyp1} and \eqref{eq:MD-basic-hyp2} we can refer to the discussion for
Example 1, with $\beta=\gamma=1$. For eq. \eqref{eq:extra-for-MD-eventually-constant} we start
noting that
\begin{multline*}
H_3(t)=\frac{1}{\sqrt{v(\delta(t))}}\left(\frac{\delta(t)D^\prime(\delta(t))}{D(\delta(t))}-\mathbb{E}[N(t)]\right)\\
=\frac{1}{\sqrt{v(\delta(t))}}\left(\sum_{k\geq 1}k\frac{d_k(\delta(t))^k}{D(\delta(t))}-
\sum_{k\geq 1}k\frac{d_{k,k}(\delta_k(t))^k}{D_k(\delta_k(t))}\left(\sum_{k\geq 0}\frac{d_{k,k}(\delta_k(t))^k}{D_k(\delta_k(t))}\right)^{-1}\right).
\end{multline*}
We remark that, by eq. \eqref{eq:initial-step} with $u=1$,
$$\sum_{k\geq 0}\frac{d_{k,k}(\delta_k(t))^k}{D_k(\delta_k(t))}=1+R_n(1,t);$$
thus we can easily check that
\begin{multline*}
H_3(t)=\frac{(1+R_n(1,t))^{-1}}{\sqrt{v(\delta(t))}}\left((1+R_n(1,t))\sum_{k\geq 1}k\frac{d_k(\delta(t))^k}{D(\delta(t))}-
\sum_{k\geq 1}k\frac{d_{k,k}(\delta_k(t))^k}{D_k(\delta_k(t))}\right)\\
=\frac{(1+R_n(1,t))^{-1}}{\sqrt{v(\delta(t))}}\left(R_n(1,t)\sum_{k\geq 1}k\frac{d_k(\delta(t))^k}{D(\delta(t))}
+\left(\sum_{k\geq 1}k\frac{d_k(\delta(t))^k}{D(\delta(t))}-
\sum_{k\geq 1}k\frac{d_{k,k}(\delta_k(t))^k}{D_k(\delta_k(t))}\right)\right)\\
=(1+R_n(1,t))^{-1}R_n(1,t)\frac{\delta(t)D^\prime(\delta(t))}{v(\delta(t))D(\delta(t))}\sqrt{v(\delta(t))}+
\frac{(1+R_n(1,t))^{-1}}{\sqrt{v(\delta(t))}}\sum_{k=1}^{n-1}k\left(\frac{d_k(\delta(t))^k}{D(\delta(t))}
-\frac{d_{k,k}(\delta_k(t))^k}{D_k(\delta_k(t))}\right).
\end{multline*}
Then we have the following statements.
\begin{itemize}
\item $\frac{d_k(\delta(t))^k}{D(\delta(t))},\frac{d_{k,k}(\delta_k(t))^k}{D_k(\delta_k(t))}\to 0$
by eqs. \eqref{eq:limit1-example-pmf-k-leq-n-1} and \eqref{eq:limit2-example-pmf-k-leq-n-1}.
\item By eq. \eqref{eq:initial-step} with $u=1$ (and by eqs. \eqref{eq:limit1-example-pmf-k-leq-n-1} 
and \eqref{eq:limit2-example-pmf-k-leq-n-1} again) 
$$R_n(1,t)=\sum_{k=0}^{n-1}\left(\frac{\lambda^k}{\Gamma(\alpha_kk+1)}\frac{(t^{\tilde{\alpha}_k})^k}{E_{\alpha_k,1}^1(\lambda t^{\tilde{\alpha}_k})}
-\frac{\lambda^k}{\Gamma(\alpha k+1)}\frac{(t^{\tilde{\alpha}})^k}{E_{\alpha,1}^1(\lambda t^{\tilde{\alpha}})}\right)\to 0;$$
actually, as it was explained for eqs. \eqref{eq:limit1-example-pmf-k-leq-n-1} 
and \eqref{eq:limit2-example-pmf-k-leq-n-1}, we can say that $R_n(1,t)\to 0$ exponentially fast, 
and therefore
$$R_n(1,t)\sqrt{v(\delta(t))}\to 0$$
because $v(\delta(t))=t^{\tilde{\alpha}/\alpha}$.
\item $\frac{\delta(t)D^\prime(\delta(t))}{v(\delta(t))D(\delta(t))}\to\Lambda^\prime(0)$ because we
can refer to the limit in eq. \eqref{eq:limit-for-mean} stated in Remark \ref{rem:limit-for-mean}
(for the previous example) with $\beta=\gamma=1$.
\end{itemize}
In conclusion $H_3(t)$ tends to zero, and therefore it is bounded. Thus eq. 
\eqref{eq:extra-for-MD-eventually-constant} is checked.

\subsubsection{A choice of the parameters for which eq. \eqref{eq:*} fails}\label{sub:counterexample}
In this section we illustrate what can happen if eq. \eqref{eq:*} fails. For simplicity we consider
the case $n=1$; however we expect to have a similar situation even if $n\geq 2$ (but the computations 
are more complicated). Thus we consider the framework in Section \ref{sub:eventually-constant-appl}
with $n=1$ and
$$\frac{\tilde{\alpha}}{\alpha}-\frac{\tilde{\alpha}_0}{\alpha_0}>0.$$
We recall that $d_0,d_{0,0}>0$. The aim is to show that, for all $\theta\in\mathbb{R}$, there exists
the limit
\begin{equation}\label{eq:counterexample-def-Psi}
\Psi(\theta):=\lim_{t\to\infty}\frac{1}{v(\delta(t))}\log\mathbb{E}\left[e^{\theta N(t)}\right]
=\lim_{t\to\infty}\frac{1}{v(\delta(t))}\log\frac{\sum_{k\geq 0}\frac{d_{k,k}(e^\theta\delta_k(t))^k}
{D_k(\delta_k(t))}}{\sum_{j\geq 0}\frac{d_{j,j}(\delta_j(t))^j}{D_j(\delta_j(t))}}\in\mathbb{R}
\end{equation}
but the function $\Psi(\cdot)$ is not differentiable and we cannot consider a straightforward application 
of the G\"artner Ellis Theorem (Theorem \ref{th:GE}), as we did in Proposition \ref{prop:LD-eventually-constant}.

Firstly we analyze $R_n(u,t)$ in eq. \eqref{eq:def-Rn}. Under our hypotheses it does not depend on $u$,
and therefore we simply write $R_1(t)$; then we have
$$R_1(t):=\frac{d_{0,0}}{D_0(\delta_0(t))}-\frac{d_0}{D(\delta(t))}
=\frac{d_{0,0}}{E_{\alpha_0,1}^1(\lambda t^{\tilde{\alpha}_0})}
-\frac{d_0}{E_{\alpha,1}^1(\lambda t^{\tilde{\alpha}})}.$$
So we can say that $R_1(t)>0$ eventually (i.e. for $t$ large enough) and $R_1(t)\to 0$ as $t\to\infty$ 
by eq. \eqref{eq:4.4.16-Gorenflo-et-al} with $u=\lambda t^{\tilde{\alpha}_0}$ and 
$u=\lambda t^{\tilde{\alpha}}$. Moreover
\begin{multline*}
\frac{1}{v(\delta(t))}\log R_1(t)=\frac{1}{t^{\tilde{\alpha}/\alpha}}
\log\left(\frac{d_{0,0}}{E_{\alpha_0,1}^1(\lambda t^{\tilde{\alpha}_0})}
-\frac{d_0}{E_{\alpha,1}^1(\lambda t^{\tilde{\alpha}})}\right)\\
=\frac{1}{t^{\tilde{\alpha}/\alpha}}
\log\left(\frac{d_0}{E_{\alpha,1}^1(\lambda t^{\tilde{\alpha}})}
\left(\frac{\frac{d_{0,0}}{E_{\alpha_0,1}^1(\lambda t^{\tilde{\alpha}_0})}}
{\frac{d_0}{E_{\alpha,1}^1(\lambda t^{\tilde{\alpha}})}}-1\right)\right)
=\frac{1}{t^{\tilde{\alpha}/\alpha}}\log\left(\frac{d_0}{E_{\alpha,1}^1(\lambda t^{\tilde{\alpha}})}\right)
+\frac{1}{t^{\tilde{\alpha}/\alpha}}\log\left(\frac{\frac{d_{0,0}}{E_{\alpha_0,1}^1(\lambda t^{\tilde{\alpha}_0})}}
{\frac{d_0}{E_{\alpha,1}^1(\lambda t^{\tilde{\alpha}})}}-1\right)
\end{multline*}
and
\begin{multline*}
\lim_{t\to\infty}\frac{1}{v(\delta(t))}\log R_1(t)\\
=\lim_{t\to\infty}\frac{1}{t^{\tilde{\alpha}/\alpha}}\log(d_0e^{-(\lambda t^{\tilde{\alpha}})^{1/\alpha}})
+\lim_{t\to\infty}\frac{1}{t^{\tilde{\alpha}/\alpha}}\log\left(\frac{d_{0,0}}{d_0}
e^{-(\lambda t^{\tilde{\alpha}_0})^{1/\alpha_0}+(\lambda t^{\tilde{\alpha}})^{1/\alpha}}-1\right)\\
=-\lambda^{1/\alpha}+\lim_{t\to\infty}\frac{-(\lambda t^{\tilde{\alpha}_0})^{1/\alpha_0}+(\lambda t^{\tilde{\alpha}})^{1/\alpha}}{t^{\tilde{\alpha}/\alpha}};
\end{multline*}
thus, by taking into account that $\frac{\tilde{\alpha}}{\alpha}-\frac{\tilde{\alpha}_0}{\alpha_0}>0$,
we get
\begin{equation}\label{eq:counterexample-lim-for-Rn}
	\lim_{t\to\infty}\frac{1}{v(\delta(t))}\log R_1(t)=0.
\end{equation}

Now we take into account eq. \eqref{eq:initial-step}.  Then, by eq. \eqref{eq:central} in 
Condition \ref{cond:eventually-constant}, for all $u>0$ we have
$$\lim_{t\to\infty}\frac{1}{v(\delta(t))}\log\frac{D(u\delta(t))}{D(\delta(t))}=\Delta(u)-\Delta(1).$$
Then we can prove the following result.

\begin{lemma}\label{lem:for-counterexample}
For all $u>0$ we have
$\lim_{t\to\infty}\frac{1}{v(\delta(t))}\log\sum_{k\geq 0}\frac{d_{k,k}(u\delta_k(t))^k}{D_k(\delta_k(t))}
=\max\left\{\Delta(u)-\Delta(1),0\right\}$.
\end{lemma}
\begin{proof}
Firstly, by eq. \eqref{eq:initial-step} with $n=1$ and by recalling that $R_1(t)>0$ 
eventually (i.e. for $t$ large enough), we can apply Lemma 1.2.15 in \cite{DemboZeitouni} and, 
by eq. \eqref{eq:counterexample-lim-for-Rn}, for all $u>0$ we have
\begin{multline*}
\limsup_{t\to\infty}\frac{1}{v(\delta(t))}\log\sum_{k\geq 0}\frac{d_{k,k}(u\delta_k(t))^k}{D_k(\delta_k(t))}\\
=\max\left\{\limsup_{t\to\infty}\frac{1}{v(\delta(t))}\log\frac{D(u\delta(t))}{D(\delta(t))},
\limsup_{t\to\infty}\frac{1}{v(\delta(t))}\log R_1(t)\right\}
=\max\left\{\Delta(u)-\Delta(1),0\right\}.
\end{multline*}
Moreover, in a similar way (actually here the application of Lemma 1.2.15 in \cite{DemboZeitouni}
is not needed), for all $u>0$ we have
\begin{multline*}
\liminf_{t\to\infty}\frac{1}{v(\delta(t))}\log\sum_{k\geq 0}\frac{d_{k,k}(u\delta_k(t))^k}{D_k(\delta_k(t))}\\
\geq\left\{\begin{array}{ll}
\liminf_{t\to\infty}\frac{1}{v(\delta(t))}\log\frac{D(u\delta(t))}{D(\delta(t))}=\Delta(u)-\Delta(1)&\ \mbox{if}\ u>1\\
\liminf_{t\to\infty}\frac{1}{v(\delta(t))}\log R_1(t)=0&\ \mbox{if}\ u\in(0,1],\\
\end{array}\right.
\end{multline*}
which yields
$$\liminf_{t\to\infty}\frac{1}{v(\delta(t))}\log\sum_{k\geq 0}\frac{d_{k,k}(u\delta_k(t))^k}{D_k(\delta_k(t))}\geq\max\left\{\Delta(u)-\Delta(1),0\right\}$$
because $\Delta(\cdot)$ is an increasing function.
\end{proof}

Finally, if we refer to the limit computed in Lemma \ref{lem:for-counterexample} 
with $u=e^\theta$, there exists the limit in eq. \eqref{eq:counterexample-def-Psi} (for all $\theta\in\mathbb{R}$) and we have
$$\Psi(\theta)=\max\left\{\Delta(e^\theta)-\Delta(1),0\right\}-\max\left\{\Delta(e^0)-\Delta(1),0\right\}
=\max\left\{\Delta(e^\theta)-\Delta(1),0\right\}.$$
Moreover, by eqs. \eqref{eq:Lambda-Lambdastar} and \eqref{eq:Lambda-for-PoganyTomovski},
we get
$$\Psi(\theta)=\max\left\{\Lambda(\theta),0\right\}=\max\left\{\lambda^{1/\alpha}(e^{\theta/\alpha}-1),0\right\}.$$
In conclusion the function $\Psi(\cdot)$ is not differentiable at the origin $\theta=0$, indeed the left 
derivative is equal to zero and the right derivative is equal to $\frac{\lambda^{1/\alpha}}{\alpha}$.

\paragraph{Acknowledgements.} The authors wish to thank the anonymous referees 
for their careful reading and suggestions to improve the presentation of the paper.
The authors also thank Roberto Garra, Roberto Garrappa, and Francesco Mainardi for
some discussion on the Prabhakar function. We also thank Camilla Feroldi for the 
activity for her thesis (which contains a preliminary version of some results in 
this paper) under the supervision of Elena Villa.

\end{document}